\documentclass[12pt]{article}

\usepackage{amssymb, latexsym, amsmath, amsfonts, amscd, fullpage, setspace}

\newtheorem{lemma}{Lemma}

\newtheorem{definition}{Definition}
\newtheorem{proposition}{Proposition}

\begin{document}

\title{A Note on Toric Varieties Associated to Moduli Spaces}
\author{James J. Uren}
\date{}

\maketitle

\begin{abstract}
In this note we give a brief review of the construction of a toric
variety $\mathcal{V}$ coming from a genus $g \geq 2$ Riemann surface
$\Sigma^g$ equipped with a trinion, or pair of pants, decomposition.
This was outlined by J. Hurtubise and L. C. Jeffrey in \cite{JH1}.
In \cite{T1} A. Tyurin used this construction on a certain
collection of trinion decomposed surfaces to produce a variety
$DM_g$ -- the so-called \emph{Delzant model of moduli space} -- for
each genus $g.$ We conclude this note with some basic facts about
the moment polytopes of the varieties $\mathcal{V}.$ In particular,
we show that the varieties $DM_g$ constructed by Tyurin, and claimed
to be smooth, are in fact singular for $g \geq 3.$

\center{MSC Primary/Secondary: 14M25/52B20}
\end{abstract}

\section{Introduction}

\subsection{Setting}

Let $\Sigma^g$ be a compact oriented two manifold of genus $g,$ and
suppose that $\partial \Sigma^g = \emptyset.$ Let
$\mathcal{M}(\Sigma^g)$ be the moduli space of gauge equivalence
classes of flat $SU(2)$ connections on $\Sigma^g.$ There is a well
known identification of $\mathcal{M}(\Sigma^g)$ with the space
$\text{Hom}(\pi_1(\Sigma^g), SU(2))/SU(2)$ of conjugacy classes of
representations of the fundamental group $\pi_1(\Sigma^g)$ into
$SU(2).$  Additionally, this space admits a symplectic form $\Omega$
(see \cite{G1} for the details.)

Now, given any simple closed loop $C$ on $\Sigma^g,$ we obtain a
function $f_C$ on the space $\text{Hom}(\pi_1(\Sigma^g),
SU(2))/SU(2),$ sending a representation class $\rho$ to

\begin{equation}
f_C(\rho) = \frac{1}{\pi}
\text{arccos}(\frac{1}{2}\text{tr}(\rho([C]))) \in [0,1]
\end{equation}

\noindent where $[C]$ is the class of $C$ in $\pi_1(\Sigma^g).$  In
\cite{JW1} Jeffrey and Weitsman proved that $f_C$ is a Hamiltonian
function for a $U(1)-$action on a large open dense subset, $U_C =
f^{-1}_C((0,1)),$ of $\text{Hom}(\pi_1(\Sigma^g), SU(2))/SU(2).$
Moreover, if $C$ and $C'$ are two simple closed loops in $\Sigma^g$
with $[C] \neq [C'],$ then the functions $f_C$ and $f_{C'}$ commute
($\{f_C, f_{C'}\}_\Omega = 0$) and the Hamiltonian flows of $f_C$
and $f_{C'}$ induce an action of a $2-$torus, $U(1) \times U(1),$ on
$U_C \cap U_{C'} \subseteq M(\Sigma^g).$ If a third loop $C''$
exists, homotopy inequivalent to $C$ and $C',$ then we obtain a
$3-$torus action, and we may continue this process provided that
additional curves can be found.

\subsection{The Case of a Trinion}

As an example, let $D$ be a trinion (a $2-$sphere with three
disjoint discs deleted) and let $C_1,$ $C_2,$ and $C_3$ be the three
boundary circles of $D.$  Denote by $\mathcal{M}(D)$ the moduli
space of gauge equivalence classes of flat $SU(2)$ connections on
$D.$

In \cite{JW1} it was shown that the map

\begin{equation}
f = (f_{C_1}, f_{C_2}, f_{C_3}): \mathcal{M}(D) \rightarrow
\mathbb{R}^3
\end{equation}
\\
\\
\noindent sends $\mathcal{M}(D)$ bijectively to the set of triples
$(x_1,x_2,x_3)$ satisfying the inequalities

\begin{align*}
x_1 + x_2 + x_3 & \leq 2 \\
x_1 + x_2 - x_3 & \geq 0 \\
x_1 - x_2 + x_3 & \geq 0 \\
-x_1 + x_2 + x_3 & \geq 0.
\end{align*}

\noindent In other words, the image of $f$ is the tetrahedron in
$[0,1]^3$ whose vertices are $(0,0,0),$ $(1,1,0),$ $(1,0,1),$ and
$(0,1,1).$ We will hereafter denote this tetrahedron by $P(D).$

\subsection{Toric Structures Coming From Trinion Decompositions}

On $\Sigma^g,$ a maximal collection of pairwise disjoint,
non-homotopic, simple closed loops has size $3g-3.$ Let $\mathcal{C}
= \{C_1, C_2, \hdots, C_{3g-3}\}$ be such a collection.  Such a set
is called a \emph{marking} of the Riemann surface $\Sigma^g.$

Recall that a \emph{trinion decomposition} of $\Sigma^g$ is a
realization of the surface as a union of $2g-2$ trinions
$\mathcal{D} = \{D_1, D_2, \hdots, D_{2g-2}\}$ glued together along
their boundary circles. Given such a decomposition, we obtain a
marking of $\Sigma^g$ by taking our set $\mathcal{C}$ to be the
collection of the $3g-3$ common boundary circles along which the
various trinions in $\mathcal{D}$ are joined.  On the other hand, it
is easy to see that any marking of $\Sigma^g$ gives rise to a
trinion decomposition of the surface.

Let us suppose that we are given a particular trinion decomposition
of $\Sigma^g.$  Let $\mathcal{D} = \{D_1, D_2, \hdots, D_{2g-2}\}$
be the set of trinions in the decomposition, and let $\mathcal{C} =
\{C_1, C_2, \hdots, C_{3g-3}\}$ be the corresponding marking of
$\Sigma^g.$  For each curve $C_i \in \mathcal{C}$ we have the
function $f_i = f_{C_i}$ (cf. equation (1)) and the set $U_i =
f^{-1}_i((0,1)).$ Let $U = \bigcap^{3g-3}_{i=1} U_i.$

We now state, without proof, two key facts from \cite{JW1} and
\cite{JW2}.

\begin{proposition}
The marking $\mathcal{C}$ determines a $3g-3$ dimensional torus $K =
\mathbb{R}^{3g-3}/\Lambda$ which acts effectively on $U,$ preserving
the symplectic form.  The lattice $\Lambda$ has rank $3g-3$ and is
spanned by the $3g-3$ standard basis vectors $e_i$ in
$\mathbb{R}^{3g-3},$ along with the vectors $g_j =
\frac{1}{2}(e_{j_1}+ e_{j_2}+ e_{j_3})$ (for $j=1, 2, \hdots,
2g-2$,) where $C_{j_1},$ $C_{j_2},$ and $C_{j_3}$ are the three
boundary circles of the trinion $D_j \in \mathcal{D}.$
\end{proposition}

\begin{proposition}
Let $f=(f_1,f_2, \hdots, f_{3g-3}): \mathcal{M}(\Sigma^g)
\rightarrow \mathbb{R}^{3g-3}.$ The restriction of $f$ to the set
$U$ is the moment map for the action of the torus $K,$ and the
closure of the image of this moment map is a convex polyhedron $P$
of dimension $3g-3.$  Let $D_j \in \mathcal{D}$ with boundary
circles $C_{j_1},$ $C_{j_2},$ and $C_{j_3},$ and denote by $\pi_j$
the projection $\mathbb{R}^{3g-3} \rightarrow \mathbb{R}^3$ defined
by $(x_1,x_2, \hdots, x_{3g-3}) \mapsto (x_{j_1},x_{j_2},x_{j_3}).$
For every $j \in \{1,2,\hdots,2g-2\}$ the image of the composition
$\pi_j \circ f$ is the tetrahedron $P(D).$  The polytope $P$ is the
intersection $\bigcap^{2g-2}_{j=1} \pi^{-1}_j(P(D)).$
\end{proposition}

Alternatively, the polytope $P$ can be described as the set of all
points $(x_1,x_2, \hdots, x_{3g-3}) \in \mathbb{R}^{3g-3}$ such that
for each $j \in \{1,2,\hdots,2g-2\}$ the triple
$(x_{j_1},x_{j_2},x_{j_3})$ satisfies the inequalities of the
previous section, where again, the indices $j_1,$ $j_2,$ and $j_3$
correspond to the three boundary circles of the $j$th trinion in
$\mathcal{D}.$  In particular, $P$ is always contained inside the
unit cube $[0,1]^{3g-3}.$

\section{From Trivalent Graphs to Toric Varieties}

Let $\Gamma$ be a trivalent graph of genus $g.$  Note that we will
allow for the possibility that $\Gamma$ has loops (edges connecting
a vertex to itself) or multi-edges (two vertices in $\Gamma$ may be
connected with more than one edge.)  Let $V(\Gamma)$ and $E(\Gamma)$
denote respectively the vertices and edges of $\Gamma.$ Counting
loops as two edges, we have $|E(\Gamma)| = 3g-3,$ and $|V(\Gamma)| =
2g-2.$

Such a graph gives us a genus $g$ $2-$manifold $\Sigma^g$ equipped
with a marking (or trinion decomposition) in the following way (see
\cite{T1}): pump up the vertices and edges of $\Gamma$ to 2-spheres
and tubes respectively.  The result is the manifold $\Sigma^g,$ and
homotopy classes of meridian circles of each of the tubes in the
pumped up graph define a set $\mathcal{C}$ of $3g-3$ disjoint,
homotopy inequivalent, simple closed loops on the surface.

Applying propositions 1 and 2 from the previous section, we obtain
from the graph $\Gamma$ a convex polytope $P(\Gamma) \in
\mathbb{R}^{3g-3},$ and a lattice $\Lambda(\Gamma)$ for the action
of the $3g-3$ dimensional torus $K(\Gamma).$  This information is,
in turn, all that is required to completely determine a toric
variety, \emph{the toric variety associated to the graph} $\Gamma,$
which we will denote by $\mathcal{V}(\Gamma).$

The toric varieties corresponding to trinion decomposed surfaces
were introduced by Jeffrey and Hurtubise in \cite{JH1}.  A primary
focus in \cite{T1} is a certain class of trivalent graphs, the
so-called \emph{multi-theta graphs}, and their corresponding toric
varieties. It is to this case that we now turn our attention.

\subsection{Multi-Theta Graphs}

The \emph{multi-theta graph of genus} $g,$ denoted $\Theta_g,$ is
best described (as in \cite{T1}) as a vertical oval $O$ crossed by
$g-1$ horizontal edges. The $2g-2$ vertices of the graph are
separated by a vertical axis of symmetry into two groups of size
$g-1.$  Each vertex is joined by an edge to the vertices immediately
above and below, and its ``twin'' opposite the axis of symmetry
(with the obvious exception of the top pair and the bottom pair,
which are connected to each other by a double edge).
\\

\noindent {\bf Example: g=2.}  In this case, our (multi-)theta graph
consists of two vertices joined by three edges.  The marking for the
corresponding surface $\Sigma^2$ consists of three curves
$\mathcal{C} = \{C_1,C_2,C_3\},$ each curve coming from some edge in
the graph.  The underlying trinion decomposition for $\Sigma^2$
consists of two trinions $D_1$ and $D_2$ glued together along their
three boundary circles.  We see that, according to proposition 2,
the three dimensional polytope $P(\Theta_2)$ is none other than the
tetrahedron $P(D)$ from the previous section.

Now, the lattice $\Lambda(\Theta_2)$ is spanned by $e_1 = (1,0,0),$
$e_2 = (0,1,0),$ and $e_3 = (0,0,1),$ together with
$g_1=g_2=(\frac{1}{2},\frac{1}{2},\frac{1}{2}).$  Let $v_i = g_1 -
e_i,$ for $i \in \{1,2,3\},$ so that $\Lambda(\Theta_2) \cong
\mathbb{Z}v_1 \oplus \mathbb{Z}v_2 \oplus \mathbb{Z}v_3.$ Evidently,
$P(\Theta_2)$ is a lattice polytope with respect to
$\Lambda(\Theta_2).$  One can verify that the $\emph{normal fan}$ of
$P(\Theta_2)$ -- the fan generated by the inward-pointing normals to
the facets of $P(\Theta_2)$ -- is a strongly convex complete
simplicial fan. More is true, for we may define an isomorphism of
the lattice $\Lambda(\Theta_2)$ with the standard lattice
$\mathbb{Z}^3$ using

\begin{equation}
A = \begin{bmatrix} 0 & 1 & 1 \\ 1 & 0 & 1 \\
1 & 1 & 0
\end{bmatrix},
\end{equation}

\noindent and additionally, $A$ maps the normal fan of $P(\Theta_2)$
to the fan generated by $(1,0,0),$ $(0,1,0),$ $(0,0,1),$ and
$(-1,-1,-1).$ This is the normal fan for the standard 3-simplex
$\Delta^3.$ It follows that the toric variety for the pair
$(P(\Theta_2),\Lambda(\Theta_2))$ is the same as the variety for the
pair $(\Delta^3, \mathbb{Z}^3),$ which is known to be
$\mathbb{C}P^3.$

\subsection{$P(\Theta_g)$ and the Variety $\mathcal{V}(\Theta_g)$}

We conclude with two straightforward facts about the polytope
$P(\Theta_g).$

\begin{lemma}
Exactly $2^g$ of the vertices of $P(\Theta_g)$ are vertices of the
unit cube $[0,1]^{3g-3}.$
\end{lemma}

\noindent \emph{Proof.}  Since $P(\Theta_g)$ is contained within
$[0,1]^{3g-3},$ if $x \in P(\Theta_g)$ and $x$ is itself a vertex of
$[0,1]^{3g-3},$ then $x$ is necessarily a vertex of $P(\Theta_g).$
Now, if $x$ is a vertex of the unit cube, then $x \in P(\Theta_g)$
if its image under each of $2g-2$ projections $\pi_j$ is a vertex of
the tetrahedron $P(D)$ (cf. proposition 2.) Such points correspond
to labellings of the edges of $\Theta_g$ with either a $0$ or a $1,$
such that for each vertex $v \in V(\Theta_g),$ the triple of edges
at $v$ are either all labelled $0,$ or exactly one is labelled $0.$

Beginning with the top pair of vertices in $\Theta_g,$ we see that
there are exactly four admissible ways to label the group of edges
emanating from the pair.  After a choice has been made for the top
pair, there are two possible labellings for the undetermined edges
adjacent to the next pair. And, for each of the remaining $g-3$
pairs of vertices there are always two possible labellings,
regardless of how the previous pair's edges were labelled.  This
gives a total of $4(2^{g-2})=2^g$ possible labellings.\\

\noindent \emph{Remark.} It must be noted that the method for
counting vertices of $P(\Theta_g)$ in the above argument is
\emph{not} exhaustive for $g \geq 3.$ That is, requiring that
$\pi_j(x)$ be a vertex of $P(D)$ for every $j \in \{1,2,\hdots,
2g-2\}$ is enough to determine that $x$ is a vertex of
$P(\Theta_g),$ but not necessary.

Proposition $5.4$ of \cite{T1} asserts that \emph{all} of the
vertices of the polytope $P(\Theta_g)$ are vertices of
$[0,1]^{3g-3},$ that there are $2^g$ in total, and that they are of
the form $(\star, \star', 0, \star, \star', 0, \star, \hdots,
\star),$ or $(\star, \star', 1, \star, \star', 1, \star, \hdots,
\star),$ where $\star$ and $\star'$ are chosen freely from
$\{0,1\}.$ We now see that this cannot be true. For example, the
above argument shows us that the point $(1,1,\hdots,1)$ cannot be a
vertex of $P(\Theta_g)$ for any $g,$ as is claimed. This can also be
seen by noting that $(1,1,\hdots,1)$ does not satisfy the
inequalities of section $1.3$ for any trinion in the decomposition of $P(\Theta_g).$ \\

In the previous section we saw that $\mathcal{V}(\Theta_2) \cong
\mathbb{C}P^3.$  One might ask whether or not any of the other
varieties $\mathcal{V}(\Theta_g)$ are also singularity free.
Unfortunately, as we shall soon see, this cannot be the case.

\begin{definition}
An $n$ dimensional polytope $P$ in $\mathbb{R}^n$ is said to be {\bf
simple} if its $1-$skeleton is an $n-$regular graph.
\end{definition}

\noindent Equivalently, an $n$ dimensional polytope $P$ is simple if
exactly $n$ facets of $P$ meet at each vertex.  It is a well-known
fact that a smooth toric variety must have a simple moment polytope.

\begin{lemma}
For every $g \geq 3,$ the polytope $P(\Theta_g)$ is non-simple.
\end{lemma}

\noindent \emph{Proof.} It follows from the proof of the previous
lemma that the origin is always a vertex of $P(\Theta_g).$  Now,
from each trinion $D_j$ in the decomposition of the underlying
surface $\Sigma^g$ we are given the set of four inequalities:

\begin{align*}
x_{j_1} + x_{j_2} + x_{j_3} & \leq 2 \\
x_{j_1} + x_{j_2} - x_{j_3} & \geq 0 \\
x_{j_1} - x_{j_2} + x_{j_3} & \geq 0 \\
-x_{j_1} + x_{j_2} + x_{j_3} & \geq 0.
\end{align*}

\noindent As we have seen, the $2g-2$ sets of inequalities of the
above type define the polytope $P(\Theta_g).$  Consider the last
three inequalities in the above set.  Each defines an affine
half-space in $\mathbb{R}^{3g-3},$ and each of these half-spaces
supports a different facet of $P(\Theta_g)$ containing the origin.
There are $3(2g-2) = 6g-6$ such facets, and so $P(\Theta_g),$ which
has dimension $3g-3,$ cannot be simple.\\

\noindent \emph{Remark.} It follows from this that the variety
$\mathcal{V}(\Theta_g)$ is singular for $g \geq 3.$  Moreover, the
above argument applies to any polytope $P(\Gamma),$ so long as
$\Gamma$ is loop-free (note that the origin is, in fact, always a
vertex of $P(\Gamma).$) So the variety $\mathcal{V}(\Gamma)$ is
singular whenever $\Gamma$ is a loop-free trivalent graph of genus
$g \geq 3.$  In proposition $5.5$ of \cite{T1} it was asserted that
not only is the polytope $P(\Theta_g)$ simple, but also that the set
of edges emanating from any vertex of $P(\Theta_g)$ forms a rational
basis for $\mathbb{R}^{3g-3}.$ In other words, it was claimed that
the polytope $P(\Theta_g)$ is \emph{Delzant,} and consequently that
the corresponding toric variety $\mathcal{V}(\Theta_g)$ -- written
there as $DM_g$ -- is always smooth. This is plainly false, since
$P(\Theta_g)$ is not even simple.

James J. Uren, \textsc{Department of Mathematics, University of
Toronto, Toronto, Ontario M5S 3G3, Canada} \\
\emph{email address:} \verb"jjuren@math.toronto.edu"

\end{document}